\documentclass[12pt]{amsart}
\usepackage{amsmath}
\usepackage{amssymb}
\usepackage{amsthm}
\usepackage{array}
\usepackage{xy}
\usepackage[pdftex]{graphicx}
\usepackage{hyperref}
\usepackage{color}
\usepackage{transparent}
\usepackage{latexsym}
\usepackage{ulem}

\usepackage{amsmath,scalefnt}
\usepackage{amsthm}
\usepackage{amssymb}
\usepackage{graphicx}
\usepackage{setspace}
\usepackage{hyperref}

\numberwithin{equation}{section}

\def\RR{\mathbb{R}}

\def\HH{\mathcal{H}}

\theoremstyle{plain}

\theoremstyle{plain}
\newtheorem{theorem}{Theorem} [section]
\newtheorem{corollary}[theorem]{Corollary}

\theoremstyle{definition}

\theoremstyle{remark}
\newtheorem{remark}[theorem]{Remark}

\numberwithin{theorem}{section}
\numberwithin{equation}{section}
\numberwithin{figure}{section}

\begin{document}
\title[]{A sharp lower bound on the polygonal isoperimetric deficit}

\author[E. Indrei]{E. Indrei $\dagger$}

 \thanks{$\dagger$\,PIRE Postdoctoral fellow}

\def\signei{\bigskip\begin{center} {\sc Emanuel Indrei\par\vspace{3mm}Center for Nonlinear Analysis\\  
Carnegie Mellon University\\
Pittsburgh, PA 15213, USA\\
email:} {\tt egi@cmu.edu}
\end{center}}

%
%\def\signln{\bigskip\begin{center} {\sc Levon Nurbekyan \par\vspace{3mm}
%Center for Mathematical Analysis,\\
%Geometry, and Dynamical Systems\\
%Departamento de Matem\'atica\\
%Instituto Superior T\'ecnico\\
%Lisboa 1049-001, Portugal\\
%email:} {\tt lnurbek@math.ist.utl.pt}
%\end{center}}

%\address{Department of Mathematics, The University of Texas, 1 University Station C1200, Austin, TX 78712}
%\keywords{ Non linear elliptic and parabolic equations, a priori estimates, maximum principle.
%\\
%\indent 2000 {\it Mathematics Subject Classification:} 35J15, 35J60, 35J62, 35J92, 35J93, 35J96, 35K10, 35K59, 35K92, 35K93, 35K96}

\makeatletter
\def\blfootnote{\xdef\@thefnmark{}\@footnotetext}
\makeatother

%\blfootnote{E. Indrei acknowledges support from the Australian Research Council, US NSF Grant DMS-0932078 administered by the Mathematical Sciences Research Institute in Berkeley, CA, and US NSF Grants OISE-0967140 (PIRE), DMS-0405343, and DMS-0635983 administered by the Center for Nonlinear Analysis at Carnegie Mellon University. L. Nurbekyan acknowledges support from the department of mathematics at the University of Texas at Austin and the Center for Mathematical Analysis, Geometry, and Dynamical Systems at Instituto Superior T\'ecnico.}
%

\date{}

\maketitle

\begin{abstract}
A sharp quantitative polygonal isoperimetric inequality is obtained.
\end{abstract}

\section{Introduction} 
%It is of interest to identify the optimal configurations in geometric and functional inequalities. Often this is carried out as a two-step process: the inequality is discovered first, and then the cases of equality is worked out. A quantitative version of a given inequality involves a ``remainder" term and simultaneously implies the inequality and exposes the minimizers. Such results, also known as ``stability estimates,"  furthermore have the advantage of quantifying in a precise way proximity of perturbations. Recent research activity has produced quantitative versions of several classical inequalities such as the isotropic isoperimetric inequality, the anisotropic isoperimetric inequality, the isoperimetric inequality inside convex cones, the Brunn-Minkowski inequality, the Sobolev inequality, and the logarithmic Sobolev inequality. The techniques involved in producing these stability results may vary greatly. For instance the main result in FMP is obtained via symmetrization techniques whereas the approach in FiMP is based on an efficient application of optimal transport theory.
The polygonal isoperimetric inequality states that if $n\ge 3$ and $P$ is an $n$-gon with area $|P|$ and perimeter $L(P)$, then the deficit is nonnegative,  
$$\delta(P):=L^2(P)-4n \tan \frac{\pi}{n}|P| \ge 0,$$
and uniquely minimized when $P$ is convex and regular. A full stability result for this classical inequality has recently been obtained in \cite{IN} via a novel approach involving a functional minimization problem on a compact manifold and the spectral theory for circulant matrices. The heart of the matter is a quantitative polygonal isoperimetric inequality for convex polygons which states that 
\begin{equation} \label{tro}
\sigma_s ^2(P) + \sigma_r^2(P) \lesssim \delta(P),
\end{equation}  
where $\sigma_s^2(P)$ is the variance of the side lengths of $P$ and $\sigma_r^2(P)$ is the variance of its radii (i.e. the distances between the vertices and their barycenter).     

The starting point of the proof is the following inequality \cite[pg. 35]{shilleto'85} which holds for any $n$-gon:
\begin{equation}\label{eq: shiletto1}
 8n^2\sin ^2 \frac{\pi}{n}\sigma_r ^2(P) \le nS(P)-4n \tan \frac{\pi}{n} |P|,
\end{equation}
where $S(P)$ is the sum of the squares of the side lengths of $P$. Since $n^2 \sigma_s ^2(P)=nS(P)-L^2(P)$, it follows that \eqref{eq: shiletto1} is equivalent to 
\begin{equation}\label{fin}
8n^2  \sin ^2 \frac{\pi}{n} \ \sigma_r ^2(P) \le   \delta(P)+n^2 \sigma_s ^2(P) .
\end{equation}
In order to establish \eqref{tro}, it is shown in \cite{IN} that  
\begin{equation} \label{do}
\sigma_s ^2(P) \lesssim \delta(P)
\end{equation}
whenever $P$ is a convex n-gon; thereafter, a general stability result is deduced via a version of the Erd{\H o}s-Nagy theorem which states that a polygon may be convexified in a finite number of ``flips" while keeping the perimeter invariant. The method of proof of \eqref{eq: shiletto1} given in \cite{shilleto'85} is based on a polygonal Fourier decomposition, whereas the technique in \cite{IN} is based on a third order Taylor expansion of the deficit (in a suitable sense) and as mentioned above involves circulant matrix theory and an optimization problem on a compact manifold. It is natural to wonder whether one can directly deduce \eqref{tro} via the method in \cite{IN} without relying on \cite{shilleto'85}. A positive answer is given in this paper. In fact, a new inequality is established which combined with \eqref{do} improves \eqref{tro}.

Let $\sigma_a^2(P)$ denote the variance of the central angles of $P$  (i.e. the angles generated by the vertices and barycenter of the vertices of $P$, see \S \ref{pre}). Then the following is true.

\begin{theorem}\label{main}
Let $n \ge 3$ and $P$ be a convex $n$-gon. There exists $c_n>0$ such that  
\begin{equation*}
c_n\, \delta(P) \ge \sigma_r^2(P)+|P|\sigma_a^2(P),
\end{equation*}
and the exponent on the deficit is sharp.
\end{theorem}

\noindent This result directly combines with \eqref{do} and yields:
 
 \begin{corollary}
 Let $n \ge 3$ and $P$ be a convex $n$-gon. There exists $c_n>0$ such that  

\begin{equation*}
c_n\, \delta(P) \ge \sigma_s^2(P)+\sigma_r^2(P)+|P|\sigma_a^2(P). 
\end{equation*}

 \end{corollary}

\begin{remark} The theorem holds for a more general class of polygons. The only requirement in the proof is that the central angles of $P$ sum to $2\pi$.    
\end{remark}

\begin{remark}
An inequality of the form $$\sigma_a^2(P) \le c_n \delta(P)$$ cannot hold in general. One can see this by a simple scaling consideration: let $P$ be a convex polygon and $P_\alpha$ be the convex polygon obtained by dilating the radii of $P$ by $\alpha>0$. Then $\delta(P_\alpha) = \alpha^2 \delta(P)$, but $\sigma_a^2(P_\alpha)=\sigma_a^2(P).$ 
\end{remark}

Quantitative polygonal isoperimetric inequalities turn out to be useful tools in geometric problems. For instance \eqref{tro} was recently utilized in \cite{CM} to improve a result of Hales which showed up in his proof of the honeycomb conjecture \cite{Hales}. Moreover, \cite{IN} has also been employed in \cite{CN} to prove a quantitative version of a Faber-Krahn inequality for the Cheeger constant of $n$-gons obtained in \cite{BF}. Related stability results for the isotropic, anisotropic, and relative isoperimetric inequalities have been obtained in \cite{FMP, FiMP, FI}, respectively.

\vskip .2in

\noindent {\bf Acknowledgements}

\vskip .05in 
\noindent The author is pleased to acknowledge support from NSF Grants OISE-0967140 (PIRE), DMS-0405343, and DMS-0635983 administered by the Center for Nonlinear Analysis at Carnegie Mellon University. Moreover, the excellent research environment provided by the Hausdorff Research Institute for Mathematics and the Rheinische Friedrich-Wilhelms-Universit\"at Bonn is kindly acknowledged.

\section{Preliminaries} \label{pre}
Let $n \geq 3$ and $P \subset \RR ^2$ be an $n$-gon with vertices $\{A_1, A_2,
\ldots, A_n \} \subset \RR ^2$ and center of mass $O$ which is taken to be the origin. For $i \in \{1,2,\ldots,n\}$, the $i$-th side length of $P$ is $l_i:=A_iA_{i+1}$, where $A_i=A_j$ if and only if $i=j$ (mod $n$); $\{r_i:=OA_i\}_{i=1}^n$ is the set of radii. Furthermore, $x_i$ denotes the angle between $\overrightarrow{OA_i}$ and $\overrightarrow{OA_{i+1}}$.

The circulant matrix method introduced in \cite{IN} is based on the idea that a large class of polygons can be viewed as points in $\mathbb{R}^{2n}$ satisfying some constraints. More precisely, consider
\begin{equation*}
  \mathcal{M}:=\Big \{ (x;r) \in \RR ^{2n}:\ x_i,r_i \geq 0,\ \eqref{eq: sum_x_i=2pi},\ \eqref{eq: sum_r^2=n},\ \eqref{eq: centroid_cond} \ \mbox{hold} \Big\},
\end{equation*}
where
\begin{equation}\label{eq: sum_x_i=2pi}
  \sum \limits _{i=1}^{n} x_i=2 \pi,
\end{equation}
\begin{equation}\label{eq: sum_r^2=n}
  \sum \limits _{i=1}^{n} r_i=n.
\end{equation}
\begin{equation}\label{eq: centroid_cond}
  \begin{cases}
  \sum \limits _{i=1}^{n} r_i \cos \left(\sum \limits_{k=1}^{i-1} x_k \right)=0,\\
  \sum \limits _{i=1}^{n} r_i \sin \left(\sum \limits_{k=1}^{i-1} x_k \right)=0.
  \end{cases}
\end{equation}
Note that $\mathcal{M}$ is a compact $2n-4$ dimensional manifold and each point $(x; r) \in \mathcal{M}$ represents a polygon centered at the origin with central angles $x$ and radii $r$; therefore, it is appropriate to name such objects \textit{polygonal manifolds}. Indeed, a point $O$ is the barycenter if and only if
\begin{equation*}
  \sum \limits _{i=1}^{n} \overrightarrow{OA_i}=0,
\end{equation*}
which is equivalent to saying that the projections of $\sum \limits _{i=1}^{n} \overrightarrow{OA_i}$ onto $\overrightarrow{OA_1}$ and $\overrightarrow{OA_1}^\perp$ vanish; in other words, $(x;r)$ satisfies \eqref{eq: centroid_cond}. Furthermore, \eqref{eq: sum_x_i=2pi} is satisfied by all convex polygons (also many nonconvex ones) and \eqref{eq: sum_r^2=n} is a convenient technical assumption which derives from scaling considerations. Note that the convex regular $n$-gon corresponds to the point $(x_*; r_*)=\left(\frac{2\pi}{n},\ldots,\frac{2\pi}{n};1,\ldots,1\right)$. 
With this in mind, the variance of the interior angles and radii of $P$ are represented, respectively, by the quantities 

$$\sigma_a^2(P)=\sigma_a^2(x;r):=\frac{1}{n} \sum \limits_{i=1}^{n} x_i^2 - \frac{1}{n^2} \left(\sum \limits_{i=1}^{n}
x_i\right)^2,$$

$$\sigma_r^2(P)=\sigma_r^2(x;r):=\frac{1}{n} \sum \limits_{i=1}^{n} r_i^2 - \frac{1}{n^2} \left(\sum \limits_{i=1}^{n}
r_i\right)^2.$$

Moreover, in $(x;r)$ coordinates, the deficit is given by the formula

\begin{align*}
\delta(P)=\delta(x; r):=\left(\sum \limits_{i=1}^{n} \left( r_{i+1}^2+r_i^2-2r_{i+1}r_i \cos
x_i\right)^{1/2}\right)^2-2n  \tan \frac{\pi}{n} \sum \limits_{i=1}^n r_i r_{i+1} \sin x_i.
\end{align*}

%
%and observe that $\mathcal{M}$ is a compact manifold of dimension $2n-4$ and all convex $n$-gons with barycenter $(0,0) \in \RR ^2$ have a representation as points $(x;r)\in \mathcal{M}$ where $(x;r)$ is associated with the $n$-gon whose $i$-th vertex has distance $r_i$ from the origin and two consecutive vertices form an angle $x_i$.

\section{Proof of Theorem \ref{main}}
By a simple reduction argument, it suffices to prove the inequality on $\mathcal{M}$: let $P$ be a convex $n$-gon and note that it is represented by $(x;r) \in \mathbb{R}^{2n},$ where $x \in \mathbb{R}^n$ denotes its interior angles and $r \in \mathbb{R}^n$ its radii. Convexity implies \eqref{eq: sum_x_i=2pi}, and \eqref{eq: centroid_cond} follows from the definition of barycenter. If $\displaystyle \sum_{i=1}^n r_i = s \neq n$, consider (by a slight abuse of notation) the polygon $P_s=(x; \frac{n}{s}r)$ obtained by scaling the radii of $P$. Evidently $\sigma_a^2(P_s) = \sigma_a^2(P),$ $|P_s| = (n/s)^2|P|$, $\sigma_r^2(P_s) = (n/s)^2 \sigma_r^2(P),$ $\delta(P_s)=(n/s)^2 \delta(P)$. Hence if the inequality stated in the theorem holds for $P_s \in \mathcal{M}$, then it also holds for $P$. 
Now let 
\begin{align*}
\phi(x; r):&= n^2 (|P|\sigma_a^2+\sigma_r^2)\\
&=\frac{1}{2}\Big(\sum \limits_{i=1}^n r_i r_{i+1} \sin x_i\Big)\Big(n\sum \limits_{i=1}^{n} x_i^2 - \big( \sum_{i=1}^n x_i \big)^2\Big)+n\sum \limits_{i=1}^{n} r_i^2 - \big( \sum_{i=1}^n r_i \big)^2,
\end{align*}
and note that it suffices to show 
\begin{equation}\label{agoal}
  \phi(x;r) \leq c \ \delta(x;r)
\end{equation}
for all $(x;r) \in \mathcal{M}.$
The polygonal isoperimetric inequality implies $\delta(x; r) \geq 0$ for every $(x; r) \in \mathcal{M}$ with $\delta(x; r) = 0$ if and only if $(x;r)=z_*:=(x_*;r_*)$. Since $\mathcal{M}$ is compact and $\delta$ is continuous it follows that 
$$
 \inf \limits_{\mathcal{M} \setminus B_{\delta}(z_*)} \delta > 0,
$$
and so \eqref{agoal} follows easily on $\mathcal{M} \setminus B_{\delta}(z_*)$. Thus it suffices to prove \eqref{agoal} for some neighborhood $B_{\delta}$ of the point $z_*$. Direct calculations imply (recall that the notation is periodic mod n) 
\begin{equation} \label{deriv}
D\phi(z_*):=(D_x \phi(z_*), D_r \phi(z_*))=0, \\
\end{equation}

\begin{equation*}
 D_{x_k x_l}\phi(z_*) =
  \begin{cases} 
      \hfill n(n-1)\sin \frac{2\pi}{n},    \hfill & k=l, \\
      \hfill -n\sin \frac{2\pi}{n}, \hfill & k\neq l, \\
  \end{cases}
\end{equation*}
\vskip .2in 
\[
 D_{r_k r_l}\phi(z_*) =
  \begin{cases} 
      \hfill 2(n-1),    \hfill & k=l, \\
      \hfill -2, \hfill & k\neq l, \\
  \end{cases}
\]
and $D_{r_kx_l} \phi(z_*)=0.$ Thus by letting $\Phi:=D^2\phi(z_*)$ it follows that
\begin{equation*}
  \Phi= \begin{pmatrix}
  n\sin \frac{2\pi}{n}\mathcal{C} & 0_{n \times n} \\
  0_{n \times n} & 2 \mathcal{C}
 \end{pmatrix},
\end{equation*}
where $0_{n \times n}$ is the $n \times n$ zero matrix and

\begin{equation*}
\mathcal{C}=\begin{pmatrix}
  n-1 &  -1 & & \cdots &  &-1\\
  -1 & n-1 &  -1 & & \cdots & \\
    & -1 & n-1 & -1 & \ddots  &   \vdots\\
   \vdots & & \ddots & \ddots & \ddots & \\
   &\vdots& \ddots & -1& n-1& -1\\
  -1 & & \cdots   &  & -1 & n-1\\
 \end{pmatrix} _{n \times n}.
 \end{equation*}
 
\noindent Moreover, $D \delta(z_*)$ is given by 
\begin{equation*}
\begin{cases}
D_{x_k}\delta(z_*)= 2n \tan \frac{\pi}{n},\\
D_{r_k}\delta(z_*)=0;
\end{cases}
\end{equation*}
hence, \eqref{eq: sum_x_i=2pi} implies  
\begin{align}
\Big \langle D\delta(z_*), (x-x_*;r-r_*) \Big \rangle&=\Big \langle D_x \delta(z_*), x-x_* \Big \rangle+\Big \langle D_r \delta(z_*), r-r_* \Big \rangle \nonumber\\
&=2n \tan \frac{\pi}{n} \sum \limits _{i=1}^{n}(x_i-(x_*)_i)=0. \label{inn}
\end{align}
 
\noindent Since $\phi(z_*)=\delta(z_*)=0$, by utilizing \eqref{deriv} and \eqref{inn} and performing a third order Taylor expansion it follows that for $z$ close enough to $z_*$, 
\begin{align} \label{low}
  \left| \phi(z)- \frac{1}{2}\langle D^2\phi(z_*)(z-z_*) , (z-z_*) \rangle \right|\leq C |z-z_*|^3,
\end{align}
and
\begin{align} \label{del}
 \left| \delta(z)- \frac{1}{2}\langle D^2 \delta(z_*)(z-z_*) , (z-z_*) \rangle \right| \leq C |z-z_*|^3, 
\end{align}
where $C>0$. In particular, there exists $\eta=\eta(n)$ such that 

\begin{equation} \label{phi}
\phi(z) \le \frac{1}{2}||\Phi||_2|z-z_*|^2+C |z-z_*|^3
\end{equation}
for all $z \in B_\eta(z_*)$. By the results of \cite[see (iv)' in \S 3 ]{IN}, it follows that
\begin{equation*}
  \inf_{w \in S_{\HH}} \langle D^2 \delta (z_*)w , w \rangle =: \sigma>0, \footnote{In fact, something stronger is proved: namely that $\displaystyle \inf_{w \in S_{\HH}} \langle D^2 f(z_*)w , w \rangle =: \sigma>0$ where $f$ is an explicit function for which $D^2f \le D^2 \delta$. This is achieved via the spectral theory for circulant matrices  and an analysis involving the tangent space of $\mathcal{M}$ at $z_*$ and the identification of a suitable coordinate system in which calculations can be performed efficiently. The barycentric condition \eqref{eq: centroid_cond} built into the definition of $\mathcal{M}$ comes up in this analysis.}
\end{equation*}
where $\HH$ is the tangent space of $\mathcal{M}$ at $z_*$ and $S_{\HH}$ is the unit sphere in $\HH$ with center $z_*$. Moreover, by continuity there exists a neighborhood $U \subset \RR ^{2n}$ of $S_{\HH}$ such that
\begin{equation*} \label{nondeg}
  \langle D^2 \delta(z_*)w , w \rangle \geq \frac{\sigma}{2},
\end{equation*}
for all $w \in U$.
Note that $\frac{z-z_*}{|z-z_*|} \in U$ for $z \in \mathcal{M}$ sufficiently close to $z_*$. Hence, there exists $\mu=\mu(\eta, \sigma) \in (0, \eta]$ such that 
$$\langle D^2 \delta(z_*)(z-z_*) , (z-z_*) \rangle \geq \frac{\sigma}{2}|z-z_*|^2$$ for $z \in B_\mu(z_*).$
In particular, for $\tilde \mu:=\min\{\mu, \frac{\sigma}{8C}\}$ and $z \in B_{\tilde \mu}(z_*),$ 
$$\delta(z) \ge \frac{1}{4}\langle D^2\delta(z_*)(z-z_*) , (z-z_*) \rangle;$$ thus, recalling \eqref{phi}, 
 $$\phi(z) \le \Big(\frac{1}{\sigma}||\Phi||_2 +\frac{2C}{\sigma} |z-z_*|\Big)\langle D^2 \delta(z_*)(z-z_*) , (z-z_*) \rangle \le c_n\delta(z),$$
where $c_n:= \frac{4}{\sigma}||\Phi||_2 +\frac{8C}{\sigma}\tilde \mu$. To achieve the second part of the theorem, it suffices to prove the existence of $c>0$ such that  
\begin{equation} \label{sha}
 \langle \Phi \, (x; r) , (x; r) \rangle \ge c|(x;r)|^2,
\end{equation}
for
\begin{equation*}
 (x;r) \in \mathcal{Z} :=\Bigg\{(x;r) :\ \sum \limits_{i=1}^n x_i=0,\ \sum \limits_{i=1}^n r_i=0 \Bigg \}.
\end{equation*}
Indeed, if \eqref{sha} holds, let $\omega:[0,\infty] \rightarrow [0,\infty]$ be any modulus of continuity (i.e. $\omega(0+)=0$) such that 
$$
\phi(z) \le c_n \omega(\delta(z)).
$$
Then for $z \in \mathcal{M}$ close to $z_*$, \eqref{del} implies 
$$\delta(z) \le c_0|z-z_*|^2,$$ for some $c_0>0$. Moreover,  $z-z_* \in \mathcal{Z}$ since $z \in \mathcal{M}$, and by combining \eqref{low} with \eqref{sha} it follows that
\begin{equation} \label{est}
\delta(z) \le c_0|z-z_*|^2 \le c_1 \langle \Phi(z-z_*), (z-z_*) \rangle  \le c_2 \phi(z) \le \tilde c \omega(\delta(z)),
\end{equation}
for some $\tilde c>0$ provided $z$ is close to $z_*$; however, since $\delta(z) \rightarrow 0$ as $z\rightarrow z_*$ and $\delta(z)>0$ for $z\neq z_*$, \eqref{est} leads to a contradiction if $$\displaystyle \liminf_{t\rightarrow 0^+} \frac{\omega(t)}{t} = 0.$$ Thus the $\displaystyle \lim \inf$ is strictly greater than zero and this implies $\omega$ is at most linear at zero.   
\noindent To verify \eqref{sha}, note first that $\mathcal{C}$ is a real, symmetric, circulant matrix generated by the vector $(n-1, -1, \ldots, -1)$. A calculation shows that the eigenvalues of $\mathcal{C}$, say $\lambda_k$, are given by 
\begin{equation} \label{eige}
\lambda_0=0 \hskip .5in \text{and} \hskip .5in \lambda _k=n \hskip .5in \text{for $k=1,\ldots, n-1$}. 
\end{equation}
Moreover, let $v_0:=(1,\ldots,1)$, and for $l \in \{1,\ldots, \lfloor \frac{n}{2} \rfloor \}$ define
   \begin{align*}
     v_{2l-1}&:=\left(1,\cos \frac{2 \pi l}{n},\cos \frac{4 \pi l}{n},\ldots,\cos \frac{2 \pi l (n-1)}{n}\right),\nonumber\\ 
     v_{2l}&:=\left(0,\sin \frac{2\pi l}{n},\sin \frac{4\pi l}{n},\ldots,\sin \frac{2 \pi l(n-1)}{n}\right).
   \end{align*} 
One can readily check that $v_k$ is an eigenvector of $\mathcal{C}$ corresponding to the eigenvalue $\lambda _{\lceil \frac{k}{2} \rceil}$, and that the set $\{v_0,v_1,\ldots,v_{n-1}\}$ forms a real orthogonal basis of $\RR ^n$ (see e.g. Proposition 2.1 in \cite{IN}). For $k=1,2,\ldots, n$, define $b_k:=(v_{k-1};0,\ldots,0) \in \RR ^{2n}$ and $b_k:=(0,\ldots,0;v_{k-n-1}) \in \RR ^{2n}$ for $k=n+1,\ldots,2n$. Since the set $\{b_k\}_{k=1}^{2n}$ forms a real orthogonal basis of $\RR ^{2n}$, given $(x; r) \in \RR ^{2n}$ there exist unique coefficients $\alpha _k \in \RR$ such that
\begin{equation*}
 (x; r)=\sum \limits _{k=1}^{2n} \alpha _k b_k.
\end{equation*}
Thus, by utilizing \eqref{eige} it follows that
\begin{align*}
  \langle \Phi(x; r),(x; r) \rangle &= \sum \limits_{k,k'=1}^{2n} \alpha _k \alpha _{k'} \langle \Phi b_k, b_{k'} \rangle\\ \nonumber
  &=n\sin \frac{2\pi}{n}\sum \limits_{k=1}^{n} \alpha _k^2 \lambda _{\lceil \frac{k-1}{2} \rceil} |b_k|^2+2\sum \limits_{k=n+1}^{2n} \alpha _k^2 \lambda _{\lceil \frac{k-n-1}{2} \rceil} |b_k|^2\\ \nonumber
  &=n^2\sin \frac{2\pi}{n} \sum \limits_{k=2}^{n} \alpha _k^2 |b_k|^2+2n\sum \limits_{k=n+2}^{2n} \alpha _k^2 |b_k|^2. \nonumber
\end{align*}  
Furthermore, if $(x; r) \in \mathcal{Z}$, $$\alpha_1=\frac{\langle (x; r), b_1\rangle}{|b_1|^2}=\sum_{i=1}^n x_i = 0,$$
$$\alpha_{n+1}=\frac{\langle (x; r), b_{n+1}\rangle}{|b_1|^2}=\sum_{i=1}^n r_i=0;$$ 
hence, 
\begin{align*}
\langle \Phi \, (x; r) , (x; r) \rangle &=n^2\sin \frac{2\pi}{n} \sum \limits_{k=1}^{n} \alpha _k^2 |b_k|^2+2n\sum \limits_{k=n+1}^{2n} \alpha _k^2 |b_k|^2\\
&\ge 2n\sum \limits_{k=1}^{2n} \alpha _k^2 |b_k|^2,
\end{align*}
and this concludes the proof.

\bibliographystyle{alpha}

\bibliography{ngonref}

\signei

\end{document}